\RequirePackage{ifpdf}
\ifpdf 
\documentclass[pdftex]{sigma}
\else
\documentclass{sigma}
\fi

\begin{document}

\allowdisplaybreaks

\renewcommand{\thefootnote}{$\star$}

\renewcommand{\PaperNumber}{020}

\FirstPageHeading

\ShortArticleName{A Bochner Theorem for Dunkl Polynomials}

\ArticleName{A Bochner Theorem for Dunkl Polynomials\footnote{This paper is a
contribution to the Special Issue ``Symmetry, Separation, Super-integrability and Special Functions~(S$^4$)''. The
full collection is available at
\href{http://www.emis.de/journals/SIGMA/S4.html}{http://www.emis.de/journals/SIGMA/S4.html}}}

\Author{Luc VINET~$^\dag$ and Alexei ZHEDANOV~$^\ddag$}

\AuthorNameForHeading{L.~Vinet and A.~Zhedanov}

\Address{$^\dag$~Centre de recherches math\'ematiques
Universite de Montr\'eal,\\
\hphantom{$^\dag$}~P.O. Box 6128,
Centre-ville Station, Montr\'eal (Qu\'ebec), H3C 3J7 Canada}
\EmailD{\href{mailto:luc.vinet@umontreal.ca}{luc.vinet@umontreal.ca}}

\Address{$^\ddag$~Donetsk Institute for Physics and Technology, Donetsk 83114,  Ukraine}
\EmailD{\href{mailto:zhedanov@fti.dn.ua}{zhedanov@fti.dn.ua}}

\ArticleDates{Received November 30, 2010, in f\/inal form February 25, 2011;  Published online February 27, 2011}

\Abstract{We establish an analogue of the Bochner theorem for f\/irst order operators of Dunkl type, that is we classify all such operators
having polynomial solutions.
Under natural conditions it is seen
that the only families of orthogonal polynomials in this category are limits of little and big $q$-Jacobi polynomials as $q=-1$.}


\Classification{33C45; 33C47; 42C05}

\begin{flushright}
\begin{minipage}{85mm}\it
This paper is dedicated to our friend\\ Willard Miller Jr
 on the occasion of his retirement
 \end{minipage}
\end{flushright}

\section{Introduction}
We have introduced recently new families of ``classical'' polynomials \cite{VZ_little,VZ_big} through limits of little and big $q$-Jacobi polynomials when $q$ goes to $-1$.  The term ``classical'' is taken to mean that the polynomials $P_n(x)$ are eigenfunctions of some dif\/ferential or dif\/ference operator $L$
\begin{gather}
L P_n(x) = \lambda_n P_n(x). \label{LPP}
\end{gather}
In the case of the novel classes of polynomials that we have identif\/ied, the operator~$L$ is of f\/irst order in the derivative operator $\partial_x$ and contains moreover the ref\/lection operator $R$ def\/ined by $R f(x)=f(-x)$. In other words the operator~$L$ can then be said to be a f\/irst order dif\/ferential-dif\/ference operator of Dunkl type~\cite{Dunkl}.

The little $-1$-Jacobi polynomials depend on~2 real parameters $\alpha$, $\beta$ such that $\alpha>-1$, $\beta>-1$. They are orthogonal on the interval $[-1,1]$ with respect to the positive weight function~\cite{VZ_little}
\begin{gather}
w(x) = (x+1)(1-x^2)^{(\alpha-1)/2}|x|^{\beta}. \label{little_w}
\end{gather}
The big $-1$-Jacobi polynomials depend on 3 real parameters $\alpha$, $\beta$, $c$ such that $\alpha>-1$, $\beta>-1$, $0<c<1$ and are orthogonal on the union of two distinct
intervals $[-1,-c]$ and $[c,1]$ with respect to the positive weight  function \cite{VZ_big}
\begin{gather}
w(x)=\frac{|x|}{x} (x+1)(x-c) (1-x^2)^{(\alpha-1)/2}\big(x^2 -c^2\big)^{(\beta-1)/2}.
\label{big_w}
\end{gather}
When $c = 0$ the weight function \eqref{big_w} becomes \eqref{little_w}. This corresponds to the well known limit relating  the big to the little $q$-Jacobi polynomials  \cite{KLS}.

Both weight functions  \eqref{big_w} and \eqref{little_w} are not associated  to ordinary classical polynomials. Ne\-vertheless, the weight function \eqref{little_w} corresponding to the little $-1$-Jacobi polynomials, belongs to the class of so-called generalized Jacobi weights (see, e.g.~\cite{Nevai}).

Like the ordinary classical orthogonal polynomials, the big and little $-1$-Jacobi polynomials possess many remarkable  properties, in particular, they can be presented in terms of Gauss' hypergeometric functions \cite{VZ_little,VZ_big}.

The purpose of the present paper is to examine in general the eigenvalue problem~\eqref{LPP}  for f\/irst-order operators $L$ of Dunkl type and hence to classify the families of polynomials that are eigenfunctions of such operators.
As a result, under the natural condition that the operator $L$ is potentially self-adjoint, we show that the big $-1$-Jacobi polynomials (and their little $-1$-Jacobi polynomials limit) exhaust the list of orthogonal polynomials with that property.

Note that these polynomials have also been found to arise from the Bannai--Ito polynomials~\cite{BI} when the orthogonality support in the latter case is extended to an inf\/inite set of points in a~limit process~\cite{VZ_big}.

\section{Dunkl-type operators of f\/irst order\\ and their polynomial solutions}

Consider the most general form of linear dif\/ferential operators~$L$ of f\/irst order which contain also the ref\/lection operator~$R$:
\begin{gather}
L= F_0(x) + F_1(x) R + G_0(x) \partial_x + G_1(x) \partial_x R,
\label{gen_L}
\end{gather}
where $F_0(x)$, $F_1(x)$, $G_0(x)$, $G_1(x)$ are arbitrary functions.

We are seeking polynomial eigensolutions of the operator~$L$, i.e.\ for every $n$ we assume that there exists
 a monic polynomial $P_n(x)=x^n + O(x^{n-1})$ which is an eigenfunction of the opera\-tor~$L$ with eigenvalue $\lambda_n$ \eqref{LPP}.
 In what follows we will assume that
\begin{gather}
  \lambda_n \ne 0  \quad \mbox{for} \ \   n=1,2,\dots,  \qquad \lambda_n \ne \lambda_m \quad \mbox{for} \ \  n \ne m. \label{ndeg_lambda}
\end{gather}
We f\/irst establish the necessary conditions for the existence of such eigensolutions.

Consider the monomials $x^n$, $n=0,1,2,\dots$. If condition \eqref{LPP} holds for all $n$, the opera\-tor~$L$ must preserve
the space of polynomials of any dimension, i.e.\ for any polynomial~$\pi(x)$ of degree~$n$, we have $L \pi(x) = \tilde \pi(x)$, where $\tilde \pi(x)$
is a polynomial of the same degree as~$\pi(x)$. Indeed, since~$\pi(x)$ can be expanded as  a f\/inite linear combination
of the polynomials~$P_i(x)$, $i=0,1,\dots, n$, the condition~$\lambda_n \ne 0$ guarantees that the polynomial $\tilde \pi(x)=L \pi(x)$
and $\pi(x)$ have the same degree.

Hence from condition \eqref{LPP} it follows that
\begin{gather}
L x^n = Q_n(x), \qquad n=0,1,2,\dots,
 \label{cond_LQ}
\end{gather} where $Q_n(x)$
are some polynomials of degree $n$:
\begin{gather*}
Q_n(x) = \lambda_n x^n + O\big(x^{n-1}\big). 
\end{gather*}
Conversely, assume that condition \eqref{cond_LQ} is valid for all monomials $x^n$, where $Q=n(x)$ are some polynomials with  leading term $\lambda_n$.

Then it is easily seen that for any $n=0,1,2,\dots$ there exists a monic polynomial $P_n(x)=x^n+O(x^n)$ satisfying the eigenvalue equation \eqref{LPP}.

Indeed, from \eqref{cond_LQ} and from conditions \eqref{ndeg_lambda} for $\lambda_n$  it follows that the operator $L$ transforms any linear space of polynomials of exact degree $n$ into itself.
Hence, from elementary linear algebra it is known that there always exists at least one eigenvector with the eigenvalue~$\lambda_n$. This eigenvector corresponds to the polynomial solution $P_n(x)$ of given degree $n$.
Moreover, again by conditions   \eqref{ndeg_lambda}, it follows that for every $n=0,1,2,\dots$ such solution is unique. We thus proved that conditions~\eqref{LPP} and~\eqref{cond_LQ} are equivalent (restrictions~\eqref{ndeg_lambda} are assumed, of course).

We shall thus proceed to verify condition~\eqref{cond_LQ} for the f\/irst values of $n$, i.e.\ $n=0,1,2,3$.
This will provide necessary conditions for the  functions $F_i(x)$, $G_i(x)$ to satisfy.

Consider the case $n=0$. Without loss of generality we can assume that
\begin{gather*}
L \{1\} =0. 
\end{gather*}
Indeed, by condition~\eqref{cond_LQ} we have $L \{1\} =c$ for some constant~$c$. Since this constant can be incorporated
additively into the function $F_0(x)$, we can assume that $c=0$. Then, from~\eqref{cond_LQ} we f\/ind that necessarily $F_1(x)=-F_0(x)$,
i.e.\ that the operator $L$ must have the form
\begin{gather}
L= F(x) (I-R) +G_0(x)\partial_x + G_1(x) \partial_x R, \label{L_F}
\end{gather}
where $F(x)=F_0(x)$ and $I$ is the identity operator. In what follows we will use the form~\eqref{L_F} for the operator~$L$.

Consider now the case $n=1$. We have
$L \{x\} = Q_1(x)$ with $Q_1(x)$ a linear function. Hence, we obtain the condition
\begin{gather}
2xF(x) + G_0(x) - G_1(x) = Q_1(x). \label{L_1_cond}
\end{gather}
Analogously, for $n=2,3$ we obtain the conditions
\begin{gather}
2x (G_0(x) + G_1(x)) = Q_2(x) \label{L_2_cond}
\end{gather}
and
\begin{gather}
2x^3F(x) + 3x^2(G_0(x)-G_1(x))=Q_3(x). \label{L_3_cond}
\end{gather}
We thus have 3 equations \eqref{L_1_cond}, \eqref{L_2_cond}, \eqref{L_3_cond} for the 3 unknown functions $F(x)$, $G_0(x)$, $G_1(x)$.
Solving these equations, we can present their general solution in the form
\begin{gather}
 G_0(x) = \frac{\mu}{x^2} + \frac{\nu_0}{x} + \rho_0 + \tau_0 x, \qquad G_1(x) = -\frac{\mu}{x^2} + \frac{\nu_1}{x} + \rho_1 + \tau_1 x, \nonumber \\
 F(x) = -\frac{\mu}{x^3} + \frac{\nu_1-\nu_0}{2x^2} + \frac{\xi}{x} + \eta, \label{sol_FG}
 \end{gather}
where $\mu$, $\nu_0$, $\nu_1$, $\rho_0$, $\rho_1$, $\tau_0$, $\tau_1$, $\xi$, $\eta$ are arbitrary constants. In what follows we will assume the restrictions
\begin{gather}
\tau_1 \ne \pm \tau_0, \qquad 2\eta + (2N+1)(\tau_0-\tau_1) \ne 0, \qquad N=0,1,2,\dots \label{res_tau}
\end{gather}
We thus have

\begin{theorem}
Let $\mu$, $\nu_0$, $\nu_1$, $\rho_0$, $\rho_1$, $\tau_0$, $\tau_1$, $\xi$, $\eta$ be arbitrary complex parameters subjected to conditions \eqref{res_tau}.
Construct the operator~$L$ using formulas~\eqref{L_F} and~\eqref{sol_FG}. Then, for any~$n$, condition~\eqref{cond_LQ} holds with
\begin{gather*}
\lambda_n = \left\{ \begin{array}{ll}
(\tau_0+\tau_1)n   \quad & \mbox{if} \   n \   \mbox{is} \  \mbox{even},\\
 2 \eta + (\tau_0-\tau_1)n   \quad  & \mbox{if} \ n \   \mbox{is} \  \mbox{odd},
\end{array}  \right . 
\end{gather*}
and  the operator $L$ has a complete set of polynomial eigenfunctions $L P_n(x) =\lambda_n P_n(x)$, $n=0,1,2,\dots$
with the nondegeneracy conditions~\eqref{ndeg_lambda}.
\end{theorem}

The proof of the theorem follows from the simple observation that for any $n=1,2,3,\dots$ we have
\[
L x^n = \lambda_n x^n + \kappa^{(1)}_n x^{n-1} + \kappa^{(2)}_n x^{n-2} + \kappa^{(3)}_n x^{n-3} = Q_n(x),
\]
where the coef\/f\/icients $\kappa^{(i)}_n$, $i=1,2,3$ are straightforwardly found from the explicit expression of the operator $L$.
Hence the operator $L$ preserves the space of polynomials. The nondegeneracy conditions \eqref{ndeg_lambda} are equivalent to the conditions~\eqref{res_tau}.

\section[Symmetrizable operators $L$]{Symmetrizable operators $\boldsymbol{L}$}

In the previous section we identif\/ied the most general form  of f\/irst-order operators $L$ of Dunkl-type that have nondegenerate polynomial
eigensolutions $P_n(x)$. Assume now that the polynomials $P_n(x)$ are orthogonal, i.e.\ that there exists a nondegenerate linear functional
$\sigma$ such that
\begin{gather}
\langle \sigma, P_n(x) P_m(x) \rangle = h_n \delta_{nm} \label{ort_cond_P}
\end{gather}
with some nonzero normalization constants~$h_n$.

This condition means that the operator $L$ is symmetric with respect to the functional $\sigma$:
\begin{gather}
\langle \sigma, L\{V(x)\} W(x) \rangle = \langle \sigma, L\{W(x)\} V(x) \rangle,
 \label{sym_L}
 \end{gather}
where $V(x)$, $W(x)$ are arbitrary polynomials. Property~\eqref{sym_L} trivially follows from~\eqref{LPP} and the completeness
of the system of polynomials $P_n(x)$.

In the case of a positive def\/inite functional, i.e.\ $h_n>0$, there is a realization of \eqref{ort_cond_P} in terms of an integral
\begin{gather*}
\langle \sigma, P_n(x) P_m(x) \rangle = \int_{a}^b P_n(x) P_m(x) d \sigma(x) 
\end{gather*}
with respect to a measure $d \sigma(x)$, where $\sigma(x)$ is a nondecreasing function of bounded variation~\cite{Ismail_book}.
The limits $a$, $b$ of integration may be f\/inite or inf\/inite.

We restrict ourselves to the case where there is a positive weight function \mbox{$w(x)\! =\! d \sigma(x)/dx \!>\!0$} inside the interval $[a,b]$
(we do not exclude the existence of concentrated masses at the end\-points~$a$,~$b$ of the interval).

Under such restrictions we have the following necessary condition for the operator~$L$:
\begin{gather}
(w(x) L)^{*} =w(x) L, \label{L_dag_L}
\end{gather}
where $M^{*}$ stands for the Lagrange adjoint operator with respect to the operator $M$.
Moreover, in what follows we will assume that all functions $F(x)$, $G_0(x)$ and $G_1(x)$ are real for the real values of the argument~$x$.
Recall that if the operator~$M$ is a dif\/ferential operator
\[
M = \sum_{k=0}^N A_k(x) \partial_x^{k}
\]
with real-valued functions $A_k(x)$ then the Lagrange adjoint operator is def\/ined by
\[
M^{*} = \sum_{k=0}^N (-1)^k \partial_x^{k} A_k(x).
\]
In the presence of the ref\/lection operator $R$, i.e.\ in the case where the operator $M$ has the form
\[
M = \sum_{k=0}^N A_k(x) \partial_x^{k} + \sum_{k=0}^N B_k(x) \partial_x^{k}R
\]
with some real-valued functions $A_k(x)$, $B_k(x)$, we def\/ine the Lagrange adjoint operator by
\begin{gather*}
M^{*} = \sum_{k=0}^N (-1)^k \partial_x^{k} A_k(x)  + \sum_{k=0}^N (-1)^k R \partial_x^{k} B_k(x) = \sum_{k=0}^N (-1)^k \partial_x^{k} A_k(x)  + \sum_{k=0}^N   \partial_x^{k} B_k(-x)R,
\end{gather*}
where we used the following formal rules
\[
R \partial_x = - \partial_x R, \qquad R^{*} = R, \qquad R B(x) = B(-x) R.
\]
These formulas assume that the interval of orthogonality is necessarily symmetric, i.e.~$b=-a$ (note that this includes possible situations
where the interval of orthogonality is the union of several intervals; in such a case, pairs of corresponding intervals should be symmetric
with respect to the ref\/lection $x \to -x$).

We thus restrict our considerations to operators $L$ which possess property~\eqref{L_dag_L} with some positive function~$w(x)$
inside the symmetric interval (or intervals) of orthogonality. Such ope\-ra\-tors are called symmetrizable while the corresponding function~$w(x)$ is called the symmetry factor for~$L$~\cite{LR,Ev_Lit}.

It is easy to see that this symmetrizability property excludes the existence of the term $G_0(x) \partial_x$ in~\eqref{gen_L} (this statement is valid only for the case of the real valued functions~$F(x)$, $G_0(x)$, $G_1(x)$  which has been assumed).
Hence we restrict ourselves in the following, to operators of the form
\begin{gather*}
L= F(x)(1-R) + G_1(x) \partial_x R  
\end{gather*}
with $G_0(x)=0$.

This implies $\mu=\nu_0=\rho_0=\tau_0=0$ and from \eqref{sol_FG}, we obtain the following expressions for the functions $G_1(x)$ and  $F(x)$:
\begin{gather*}
G_1(x) = \frac{\nu_1}{x} + \rho_1 + \tau_1 x, \qquad F(x) = \frac{\nu_1}{2x^2} + \frac{\xi}{x} + \eta 
\end{gather*}
with 5 arbitrary parameters $\nu_1$, $\rho_1$, $\tau_1$, $\xi$, $\eta$.

Note that we can multiply the operator~$L$ by an arbitrary nonzero constant $L \to \kappa_0 L$.
A~sca\-ling transformation of the independent variable~$x \to \kappa_1 x$ with an arbitrary nonzero constant~$\kappa_1$ can also be performed.
Using this freedom, without loss of generality,  we can always reduce the number of arbitrary parameters to~3.

Assume that $\nu_1 \tau_1 \ne 0$. We can then rewrite $G_1(x)$ in the form
\begin{gather}
G_1(x)= \frac{g_1(x-d)(x+c)}{x},
\label{G_1_quad}
\end{gather}
where $g_1$, $d$, $c$ are real constants. Under the condition that $d \ne c$, we can set $g_1=2$ and $d=1$ by an appropriate choice
of the parameters $\kappa_0$, $\kappa_1$
so as to  have
\begin{gather}
G_1(x) =\frac{2(x-1)(x+c)}{x} \label{G_1_big}
\end{gather}
with one arbitrary real parameter $c$. We will assume that $0<c<1$. (The cases $c>1$ or $c<-1$ can be analyzed analogously.)
For the function $F(x)$ we have{\samepage
\begin{gather}
F(x) = -\frac{c}{x^2} + \frac{\beta-\alpha c}{x} -\alpha-\beta-1 \label{bJ_param}
\end{gather}
with two arbitrary  real parameters $\alpha$, $\beta$, $c$.}

In \cite{VZ_big} it was shown that the operator $L$ with the coef\/f\/icients $G_1(x)$, $F(x)$ given by~\eqref{G_1_quad} and~\eqref{bJ_param} have big $-1$-Jacobi polynomials as eigensolutions. In the limiting case $c=0$ the little $-1$-Jacobi polynomials are obtained.

Hence we have
\begin{theorem}
The functions $G_1(x)$ and $F(x)$ given by  \eqref{G_1_quad} and \eqref{bJ_param} thus provide the most general first order operator of
Dunkl-type under our hypotheses. They define the operator $L$ that has the big $-1$-Jacobi polynomials as eigenfunctions
\end{theorem}

Note that we have considered the generic (nondegenerate) choice of real parameters in the functions~$G_1(x)$ and $F(x)$.
In the next section we will examine the  symmetrizability property of the operator $L$ and consider also the degenerate cases
with respect to the parameters.

\section{Pearson-type equation for the weight function}

It is instructive to derive directly the expression for the weight function $w(x)$ from condition~\eqref{L_dag_L}.

Denote $M=w(x)L=w(x)F(x)(1-R) +w(x) G_1(x)\partial_x R$. Assuming that all parameters of the functions $F(x)$, $G_1(x)$, $w(x)$ are real, we have
\begin{gather*}
M^{*} = w(x) F(x) - w(-x)F(-x) R - \partial_x R w(x) G_1(x)\\
\phantom{M^{*}} {}=
w(x) F(x) - w(-x)F(-x) R +  \frac{d(w(-x) G_1(-x))}{d x}R +
w(-x) G_1(-x)\partial_x R.
\end{gather*}
Condition $M^{*}=M$ thus means that the following two equations hold:
\begin{gather}
w(x)G_1(x)= w(-x) G_1(-x) \label{wG}
\end{gather}
and
\begin{gather}
w(-x)F(-x)-w(x)F(x) = \frac{d}{d x} w(-x) G_1(-x) =  \frac{d}{d x} w(x) G_1(x).\label{wF}
\end{gather}
With $G_1(x)$ given by \eqref{G_1_big}, we have from \eqref{wG}
\begin{gather}
w(x) = \theta(x) (x+1)(x-c) W\big(x^2\big), \label{wW}
\end{gather}
where $\theta(x) = x/|x|$ is the standard sign function and $W(x^2)$ is an arbitrary even function of~$x$.
We will suppose that $W(x^2)>0$ inside the interval $[c,1]$. Assume that $\alpha>-1$, $\beta>-1$.
Then from~\eqref{wW} it follows that the function~$w(x)$ is positive inside the two symmetric intervals $[-1,-c]$ and $[c,1]$.

Substituting expression~\eqref{wW} into the second condition of~\eqref{wF} and using the explicit expression~\eqref{bJ_param},
we obtain a dif\/ferential equation for the function~$W(y)$:
\[
2(y-1)\big(y-c^2\big)W'(y) + \big((2-\alpha-\beta)y + \beta-1 + c^2(\alpha-1)  \big)W(y) =0,
\]
which has the following general solution
\[
W(y) = C (1-y)^{(\alpha-1)/2} \big(y-c^2\big)^{(\beta-1)/2}
\]
with an arbitrary constant $C$. Thus the function $w(x)$ has the expression
\begin{gather*}
w(x) = C \theta(x) (x+1)(x-c) (1-x^2)^{(\alpha-1)/2} \big(x^2-c^2\big)^{(\beta-1)/2}. 
\end{gather*}
Assuming $C>0$ and $\alpha>-1$, $\beta>-1$, we see that $w(x)>0$ inside the two  symmetric intervals $[-1,-c]$ and $[c,1]$.
This coincides with the weight function~$w(x)$~\eqref{big_w} corresponding to the big $-1$-Jacobi polynomials which was derived in~\cite{VZ_big}
from a completely dif\/ferent approach.

Consider now possible degenerate cases of the function $G_1(x)$.

$(i)$ Assume that $\nu_1=0$ but $\rho_1 \tau_1 \ne 0$. Then using scaling transformations, we can always reduce the functions $G_1(x)$ and $F(x)$ to
\[
G_1(x)=2(1-x), \qquad F(x) = \alpha+\beta+1-\beta/x.
\]
From \eqref{wG} we get $w(x) = (x+1) W(x^2)$ with some function~$W(y)$. Then from \eqref{wW} we obtain a~Pearson-type dif\/ferential
equation or the function $W(y)$:
\begin{gather*}
2y(y-1) W'(y) + ((1-\alpha-\beta)y +\beta)W(y) =0. 
\end{gather*}
Its general solution is
\begin{gather*}
W(y) = C y^{\beta/2}(1-y)^{(\alpha-1)/2}. 
\end{gather*}
Hence we obtain for the weight function
\begin{gather*}
w(x) = C (x+1) |x|^{\beta} \big(1-x^2\big)^{(\alpha-1)/2}. 
\end{gather*}
This weight function
corresponds to the one of the little $-1$-Jacobi polynomials  \eqref{little_w} which
are orthogonal on the interval $[-1,1]$~\cite{VZ_little}.
Equivalently, this case corresponds to setting $c=0$ in formulas~\eqref{G_1_big} and~\eqref{bJ_param}.

$(ii)$ Assume that $\nu_1=\rho_1=0$ but $\tau_1 \ne 0$. Then we can write without loss of generality
\[
G_1(x)=2x, \qquad F(x) = \alpha+\beta+1-\beta/x.
\]
From \eqref{wG} we now have $w(x) = \theta(x) W(x^2)$ and from \eqref{wF} follows the Pearson equation
\[
2y W'(y) + (\alpha+\beta+2)W(y)=0
\]
with solution $W(y)=C y^{-(\alpha+\beta+2)/2}$. Whence
\begin{gather}
w(x) = C \theta(x) |x|^{-(\alpha+\beta+2)/4}. \label{w_3}
\end{gather}
The function $w(x)$ of~\eqref{w_3} is not positive however, in any symmetric interval $[-a,a]$ (or any union of symmetric intervals).
Hence, this case does not lead to orthogonal polynomials with positive measure on the real axis.

$(iii)$ Consider the degenerate case when the numerator in expression \eqref{G_1_quad} has coinciding zeros. Using scaling transformations,
$G_1(x)$ and $F(x)$ can be cast in the form
\begin{gather*}
G_1(x)= \frac{2(x-1)^2}{x}, \qquad F(x) = x^{-2} + ax^{-1} +b 
\end{gather*}
with arbitrary constants $a$, $b$.

From condition \eqref{wG} we f\/ind
\begin{gather*}
w(x) = \theta(x) (x+1)^2 W\big(x^2\big) 
\end{gather*}
with an arbitrary even function $W(x^2)$. From condition \eqref{wF} we then obtain a dif\/ferential equation for the function $W(y)$:
\[
2(y-1)^2 W'(y) + \left((b+3)y + 2a+b-1 \right)W(y)=0
\]
with general solution
\begin{gather*}
W(y) = C (y-1)^{-(b+3)/2}   \exp\left( \frac{a+b+1}{y-1} \right) 
\end{gather*}
with an arbitrary constant $C$.
Thus the weight function can be presented in the form
\begin{gather*}
w(x)= C \theta(x) (x-1)^2 \left(x^2-1\right)^{-(b+3)/2} \exp\left( \frac{a+b+1}{x^2-1} \right). 
\end{gather*}
It is  easily seen that due to presence of the sign function $\theta(x)$ there are no symmetric intervals $[-d,d]$
(or pairs of symmetric intervals) inside which the function $w(x)$ can be positive.
Hence in this (degenerate) case the operator $L$ is not symmetrizable.

$(iv)$ Consider the case $\tau_1=0$, $\nu_1 \rho_1 \ne 0$. We can here posit $G_1(x) = 2(1-x^{-1})$ and $F(x)=-x^{-2} - \alpha x^{-1} -\beta-1$.
This leads to $w(x) = \theta(x) (x+1) W(x^2)$ and for the function $W(y)$ we have the equation
\[
2(1-y)W'(y) + (\alpha+\beta) W(y) =0,
\]
whence
\[
w(x) = C \theta(x) (x+1) \big(1-x^2\big)^{(\alpha+\beta)/2}.
\]
Again it is impossible to realize the condition $w(x)>0$ on a symmetric interval.

$(v)$ Finally consider the case when $\tau_1=\rho_1=0$. Then $G_1(x)
= -2/x$, $F(x) = -x^{-2} -\alpha x^{-1} - \beta$. This case
corresponds to the weight function
\[
w(x) = C \theta(x) \exp\left(-\frac{\beta}{2}x^2 \right)
\]
and here also it is impossible to have a positive weight function on a symmetric interval.

We see that in all cases, the function $w(x)$ possesses a factorization $w(x)=\pi(x) W(x^2)$, where $\pi(x)$ is a polynomial of f\/irst or second degree,
and the function $W(y)$ obeys a Pearson-type equation that we have derived.
Roughly speaking, the function $W(y)$ resembles the corresponding weight function for the classical orthogonal polynomials~\cite{NSU}.
But, in contrast to the classical case, no analogs of Laguerre or Hermite polynomials appear in the case $q=-1$.

\section{Conclusions}

We showed that the symmetrizability property of the Dunkl operator $L$ with a positive weight function $w(x)$
 on a symmetric interval (or a union of symmetric intervals) leads to very strong restrictions upon the associated ``classical'' orthogonal polynomials.
 Namely, only  big and little $-1$-Jacobi polynomials admit Dunkl analogs. A number of generalizations of\/fer themselves.

For example it was shown in \cite{Cheikh} that the generalized Hermite polynomials $H_n^{(\mu)}(x)$ and the generalized Gegenbauer polynomials $S_n^{(\xi, \eta)}(x)$ (see, e.g.~\cite{Chi,Ros, Belm} for their def\/inition and pro\-perties) are eigenfunctions
\[
L_{\mu} P_n((x) = \lambda_n P_n(x)
\]
of an operator $L_{\mu}$ which is quadratic
\begin{gather}
L_{\mu} = \sigma(x) T_{\mu}^2 + \tau(x) T_{\mu}
\label{L_T}
\end{gather}
with respect to the Dunkl operator
\[
T_{\mu}=\partial_x + \mu x^{-1} (I-R).
\]
Here $\sigma(x)$ is a polynomial of degree not exceeding 2 and $\deg(\tau(x))=1$.

Hence the generalized Hermite and Gegenbauer polynomials can be considered as ``classical'', but in contradistinction to the general situation considered  here, they satisfy an eigenvalue equation of second order in the Dunkl operator. The classif\/ication of all second-order Dunkl-type operators having orthogonal polynomials as solutions is a more involved problem.

Note that the operator $L$ def\/ined by \eqref{L_T} is a direct generalization of the classical hypergeometric operators{\samepage
\[
L = \sigma(x) \partial_x^2 + \tau(x) \partial_x
\]
leading to the classical orthogonal polynomials as eigensolutions~\cite{NSU}.}

Let us mention also that  operators of Dunkl-type have also been used in \cite{KB} and \cite{Chouchene}
 to construct some polynomial eigenvalue solutions. In these cases however,  the Dunkl operators are not
 symmetrizable, and the solutions therefore do not belong to the class of ordinary orthogonal polynomials.

\subsection*{Acknowledgments}

 The authors are indebted to R.~Askey, C.~Dunkl, T.~Koornwinder, W.~Miller, V.~Spiridonov, P.~Terwilliger and P.~Winternitz for stimulating communications and to an anonymous referee for valuable remarks. AZ thanks CRM (U de Montr\'eal) for its hospitality.

\pdfbookmark[1]{References}{ref}
\LastPageEnding

\end{document}